\documentclass[12pt,letterpaper]{amsart}
\usepackage{amsmath,amsthm,amsfonts,amssymb,amscd}
\usepackage{lastpage}
\usepackage{enumerate}
\usepackage{fancyhdr}
\usepackage{mathrsfs}
\usepackage{xcolor}
\usepackage{graphicx}
\usepackage{listings}
\usepackage{hyperref}
\usepackage{enumitem}
\usepackage{tikz}
\usepackage{comment}

\usepackage[colorinlistoftodos]{todonotes}

\usepackage{parskip}
\newtheorem{thm}{Theorem}[section]

\newtheorem*{introthm*}{Theorem}

\newtheorem*{lem*}{Lemma}

\theoremstyle{definition}

\theoremstyle{remark}
\newtheorem{rem}[thm]{Remark}
\newtheorem*{notation*}{Notation}

\newcommand{\F}{{\mathbb F}}

\renewcommand{\P}{{\mathbb P}}

\newcommand{\Q}{{\mathbb Q}}
\newcommand{\Z}{{\mathbb Z}}
\newcommand{\C}{{\mathbb C}}

\newcommand{\CP}{{\mathbb C \mathbb P}}

\newcommand{\cO}{{\mathcal O}}

\author[L. Borisov]{Lev Borisov}
\address{Hill Center
Department of Mathematics, Rutgers University, NJ 08854}
\email{\href{mailto:borisov@rutgers.edu}{borisov@rutgers.edu}}

\author[Z. Lihn]{Zachary Lihn}
\address{Department of Mathematics, Columbia University, New York, NY 10027, USA}
\email{\href{mailto:zal2111@columbia.edu}{zal2111@columbia.edu}}

\title[Fake Projective Plane in $\P^5$]{Realizing a Fake Projective Plane as a Degree 25 Surface in $\P^5$}

\begin{document}
\maketitle
\begin{abstract}
    Fake projective planes are smooth complex surfaces of general type with Betti numbers equal to that of the usual projective plane. Recent explicit constructions of fake projective planes embed them via their bicanonical embedding in $\mathbb P^9$. In this paper, we study Keum's fake projective plane $(a=7, p=2, \{7\}, D_3 2_7)$ and use the equations of \cite{Borisov} to construct an embedding of fake projective plane in $\mathbb P^5$. We also simplify the 84 cubic equations defining the fake projective plane in $\mathbb P^9$.
\end{abstract}
\section{Introduction}
The Enriques–Kodaira classification splits compact complex surfaces $S$ into 10 classes based largely on their Kodaira dimension $k(S)$. While surfaces with Kodaira dimension $<2$ are better understood, those of general type with maximum Kodaira dimension $k(S)=2$ still need a detailed classification.

To each minimal model of a surface $S$ one associates a triple of numerical invariants $(p_g,q,K_S^2)$, where $p_g=h^0(S,K_S)$ is the geometric genus, $q=h^1(S, \cO_S)$ is the irregularity, and $K_S^2$ is the self-intersection number of the canonical class $K_S$. These determine all the other classical invariants such as the topological Euler characteristic $e_{top}(S)=12\chi(\cO_S)-K^2_S$ and the plurigenera $P_m(S)=h^0(S,mK_S)$ \cite{Hartshorne}. It turns out that producing surfaces with low $p_g$ and $q$ is quite difficult and a complete classification appears far away \cite{Borisov and Fatighenti}. In the case of $p_g=q=0$, one has the Bogomolov-Miyaoka-Yau inequality  $K^2_S\leq 9$. The focus of this paper is the extreme case of surfaces with $p_g=q=0$ and $K^2_S=9$. These are the \textit{fake projective planes} (often called FPPs for short) which by definition are complex projective surfaces of general type with Hodge diamond 
	$$
	\begin{array}{ccccc}
	&&1&&\\
	&0 &&0&\\
	0&&1&&0\\
	&0&&0&\\
	&&1&&
	\end{array}
	$$
which is the same as that of $\CP^2$. The existence of a fake projective plane was first proved by Mumford \cite{Mumford} by expressing the surface as a quotient of a 2-adic analog of the complex two-dimensional ball 
\[
\mathcal B ^2 = \{(z_1, z_2) \in \C^2: |z_1|^2 + |z_2|^2 \leq 1\}
\]
by a finitely generated group.

The general theory ensures that each fake projective plane is algebraic. By Noether's formula we know that $c_1^2=9$ and so all FPPs have $c_1^2=3c_2=9$, where $c_1, c_2$ are the Chern numbers. This implies that each FPP is a quotient of $\mathcal B ^2$ by an infinite discrete group \cite{Yau}. These ball quotients are determined by their fundamental group up to holomorphic or anti-holomorphic isomorphism \cite{Mostow} and come in complex conjugate pairs \cite{Kharlamov and Kulikov}. Each of the groups are arithmetic \cite{Klinger} and come in a finite list of classes \cite{Prasad and Yeung}. 

Based on the work of Prasad and Yeung  \cite{Prasad and Yeung}, a complete classification was obtained by Cartwright and Steger \cite{Cartwright and Steger}. All fake projective planes are quotients of $\mathcal B^2$ by explicit co-compact torsion-free arithmetic subgroups of $\text{PU} (2,1)$. The classification was accomplished with significant use of computer calculations. There are 50 conjugate pairs of fake projective planes split among 28 classes. Each FPP is a ball quotient $\mathcal B^2 / \Gamma$ where $\Gamma$ is the fundamental group, and where the automorphism group is $ N(\Gamma)/\Gamma $ with $N(\Gamma)$ the normalizer of $\Gamma$ in $\text{PU} (2,1)$. The torsion of the Picard group of $\P^2_{fake}$ is equal to the abelianization of $\Gamma$. Various cover relations between related surfaces are also known \cite{Cartwright and Steger}.

\subsection{The Geometry of Keum's Fake Projective Plane}
In this paper, we will focus on the fake projective plane $(a=7, p=2, \{7\}, D_3 2_7)$ in Cartwright-Steger classification. First constructed in \cite{Keum 6}, it is named Keum's fake projective plane and we will denote it by $\P^2_{Keum}$. Its automorphism group has maximum order among all FPPs, being equal to the semi-direct product of a normal cyclic subgroup $C_7$ of order 7 and a non-normal cyclic subgroup $C_3$ of order 3. By the Cartwright-Steger classification, there are three other fake projective planes in its class including Mumford's first fake projective plane.

For the rest of this paper, we will let $K$ denote the canonical class of Keum's fake projective plane. The minimal resolution $Y$ of the quotient $\P^2_{Keum}/C_7$ by the subgroup $C_7$ of its automorphism group has interesting geometry which we describe briefly. 

Recall that a singular point of type $\frac{1}{m}(1,a)$ is a cyclic quotient singularity given local analytically by the action $(x,y) \mapsto(\zeta x, \zeta^a y)$ on $\C^2$ for $\zeta$ a primitive $m$th root of unity. $Y$ has three singular points of type $\frac{1}{7}(1,3)$ 
permuted by the residual $C_3$ automorphism group of $\P^2_{Keum}$. It is also a Dolgachev surface fibered over $\P^1$, with generic fibers of genus 1, two multiple fibers, three nodal fibers, and one fiber of type $I_9$. The two multiple fibers are $2F_3$ and $3F_2$, which have multiplicity 2 and 3 respectively. The reductions $F_3$ and $F_2$ are linearly equivalent to $3K_Y$ and $2K_Y$. We refer to \cite{Keum 6, Borisov} for more details.

\subsection{Explicit Construction of \texorpdfstring{$\P^2_{Keum}$}{P2Keum}}
In \cite{Borisov}, Keum's fake projective plane was explicitly constructed via its bicanonical embedding as the vanishing set of 84 cubic equations in $\P^9$. One first constructs a birational model $Y_0$ of $Y$ as a system of quadrics in 8 variables defined over $\Q(\sqrt{-7})$. Included is a construction of the double and triple fibers and the $C_3$ action on $Y_0$. A degree 7 extension of the field of rational functions of $Y_0$ gives the sevenfold cover of $Y_0$, which is exactly $\P^2_{Keum}$. Ten sections of $\cO(2K)$ are constructed from this description and the embedding in $\P^9$ is finally given by 84 cubic equations in the 10 variables $P_0, \ldots, P_9$. 

A perennial question is how to simplify the equations of a fake projective plane, which can have polynomials with coefficients hundreds to thousands of decimal digits long. In this paper, we give a simplified version of the equations of Keum's fake projective plane in \cite{Borisov}. We use the equations to find an embedding of $\P^2_{Keum}$ as a degree 25 surface in $\P^5$. The embedding is given by sections of $\cO(5H)$, where $H$ is a divisor such that $3H$ is linearly equivalent to $K$. Finally, we exhibit the surface as a system of 56 sextics in $\P^5$ with coefficients in the field $\Q(\sqrt{-7})$. 

The paper is organized as follows. In Section \ref{section2} we outline the steps to simplify the 84 cubics defining $\P^2_{Keum}$ in $\P^9$. We follow the strategy described in \cite{Borisov} by explicitly calculating the nonreduced linear cuts on $\P^2_{Keum}$ corresponding to $2$-torsion in the Picard group. Using these equations, in Section \ref{section3} we describe the steps to embed $\P^2_{Keum}$ in $\P^5$. Specifically, we compute global sections of $\cO(5H)$ as global sections of the divisor $18H-9H-4H$ and explain the key idea that allowed us to find $H^0\left(\P^2_{fake}, \cO(4H)\right)$. Section \ref{section4} concludes with future directions.

\begin{rem}
A defining feature of recent constructions of fake projective planes is their heavy use of computer algebra software. To that end, this project depended heavily on the use of the Mathematica software system \cite{Math} and the computer algebra systems Magma \cite{Magma} and Macaulay2 \cite{Macaulay2}.
\end{rem}

\begin{rem}
The 84 cubics in $\P^9$ and the 56 sextics in $\P^5$ are still too large to be included in the printed paper.
\end{rem}

\section{Simplification of Keum's Fake Projective Plane} \label{section2}
We will begin by simplifying the explicit equations of Keum's fake projective plane $\P^2_{Keum}$ found in \cite{Borisov}. This is done by looking for nonreduced cuves on $\P^2_{Keum}$ which correspond to $2$-torsion in the Picard group. We proceed by making a coordinate change that makes the curve nicer in our new basis.

\subsection*{Step 1: Finite Field Search for Nonreduced Curves}

By the Cartwright-Steger classification, the torsion in the Picard group of $\P^2_{Keum}$ is $C_2^3$. In addition, the automorphism group is $ C_7 \rtimes C_3$, the semidirect product of $C_7$ and $C_3$.


We claim that 2-torsion classes give nonreduced curves in $|2K|$. Let $L$ be a 2-torsion class in the Picard group. By \cite{GKS}, we have $h^0(\P^2_{Keum}, K+L) = 1$. Hence, up to scaling, there is a unique section $s_{L}\in H^0(\P^2_{Keum}, K+L$). The square of $s_L$ is in $H^0(\P^2_{Keum}, 2K)$ and gives rise to a nonreduced curve. 

We will further assume that the nonreduced curve is $C_3$ invariant. This reduces the search to nonreduced curves of the form 
\[
a_0 P_0 + a_1(P_1+P_2+P_3) + a_2(P_4+P_5+P_6)+a_3(P_7+P_8+P_9)
\]
up to scaling (so we subsequently set $a_0=1)$. 
To look for such curves we look at a finite field reduction of $\P^2_{Keum}$ over $\F_p$ for suitable $p$. More precisely, such suitable $p$ contains a square root of $-7$ and has the same Hilbert polynomial for $\P^2_{Keum}$ modulo $p$. We picked $p=43$ with $\sqrt{-7} \equiv 6 \mod{43}$ which was an arbitrary small prime with the aforementioned conditions. Using Magma, we ran an exhaustive search for all $a_1, a_2, a_3$ in $\F_{43}$ and checked if the corresponding curve is nonreduced. We obtained the curve 
\[
P_0 + 24 (P_1 + P_2 + P_3) + 0 (P_4 + P_5 + P_6) + 28 (P_7 + P_8 + P_9).
\]

\subsection*{Step 2: Lift to Characteristic 0}
We lift this curve to $\Q(\sqrt{-7})$ as follows. Using Magma, we calculate some points in $\F_{43}$ lying on $\P^2_{Keum}$ and the nonreduced curve. We then apply a variant of Hensel lifting to lift the curve to $\Z/43^k\Z$ for higher $k$ at each step, obtaining a $p$-adic approximation. 

The lifting process was done by finding, at each point, two linearly-independent tangent vectors in $\P^9(\F_{43})$ that are orthogonal to all polynomials defining $\P^2_{Keum}$ and the linear cut. We modified the points, tangent vectors, and the linear cut at each stage to lift them to higher powers of $43$ such that the orthogonality conditions held; this reduced to solving a system of linear equations modulo $43$. After a sufficiently high power of $43$ we identify the corresponding algebraic numbers by applying a lattice reduction algorithm. We obtain the curve 
\begin{align*}
&P_0 + \frac{(-1 + \sqrt{-7})}{2}  (P_1 + P_2 + P_3) + \frac{(272 - 848 \sqrt{-7})}{7}
  (P_4 + P_5 + P_6) \\
  &+
   \frac{(832 - 192 \sqrt{-7}) }{7}(P_7 + P_8 + P_9)
\end{align*}
which we verify is nonreduced numerically. 

Thus we have found one nontrivial $C_3$-invariant torsion line bundle. It is not $C_7$-invariant because the corresponding nonreduced linear cut is not $C_7$ invariant. Its orbit therefore has $7$ elements, which combined with knowledge of the torsion of the Picard group as $ C_2^3$ shows that the action of the automorphism group on the torsion in Picard group is transitive.

\subsection*{Step 3: Setting Up the Coordinate Change}
Finally we set up the coordinate change to find a nicer basis for $H^0(\P^2_{Keum}, 2K)$ in order to simplify the equations defining the fake projective plane. We use a coordinate change from $P_i$ to $Q_i$ that respects the automorphisms on the surface such that the the nonreduced cut becomes 
\[
Q_0 + Q_1+Q_2+Q_3+Q_7+Q_8+Q_9.
\]

These conditions leave one free parameter in the coordinate change. We fix the free parameter by choosing it in such a way to set the "simplest" coefficient in the equations to 1. This allows us to find a version of the 84 equations with significantly smaller coefficients. 

We simplify the equations further by reducing the number of monomials in the equations. We take random linear combinations of the seven equations in each $C_7$ weight and select those that span the space and have the fewest monomials.

\section{Embedding of a fake projective plane into \texorpdfstring{$\P^5$}{P5}} \label{section3}
In this section, we will describe the process that led us to find the equations of an embedding of $\P^2_{Keum}$ in $\P^5$. 
Let $H$ be a divisor such that $3H = K$, where $K=K_{\P^2_{Keum}}$ is the canonical divisor of $\P^2_{Keum}$. Calculations of $h^0(\P^2_{Keum}, nH)$  show that $|5H|$ has the expected dimension such that the corresponding map to projective space is $\P^5$. Thus we aim to construct $|5H|$ explicitly, which will give the desired map $\P^2_{Keum} \to \P^5$.

\begin{table}[!h]
\begin{center}
\begin{tabular}{|c||c|c|c|c|c|c|c|c|c|c|c|c|}
\hline 
   $ n$ & 3&4&5&6&7&8&9&10&11&12 \\
   \hline
    $h^0\left(\P^2_{Keum},nH\right)$ &0&3&6&10&15&21&28&36&45&55  \\
    \hline
\end{tabular}

\caption{Dimensions of $H^0(\P^2_{Keum}, nH)$ for different values of $n$, where $3H=K$}.
\end{center}
\end{table}

Recall that $Y$ denotes the quotient $\P^2_{Keum}/C_7$ of Keum's fake projective plane by its $C_7$ automorphism subgroup. It has residual automorphism group $C_3$ and has a double fiber $F = 3K_Y$.

We will construct $|5H|$ as the space $|18H-9H-4H|$. We first find $|9H|=|18H-9H|$ by expressing 112 cubic equations in the $Q_i$ (which lie in $18H=6K$) that vanish on $9H$. Crucial to this construction is the preimage of the double fiber $F$ of $Y$ which we use to find points on $9H$. We then compute $|4H|$. This required the use of several important ideas which are detailed in Step 2 below. Finally, after constructing $4H$ we may find $5H$ as linear combinations of the equations of $9H$ vanishing on $4H$. We conclude by using the explicit equations in $\P^5$ to reconstruct the $C_3$ action on $\P^2_{Keum}$ in its embedding into $\P^5$.

\subsection*{Step 1: Constructing $|9H|$}

The preimage of the double fiber on $Y$ has divisor class $3K = 9H$ \cite{Borisov}. Hence to construct $|9H|$ we are led to find polynomials on $Y$ vanishing on the double fiber. Recall that \cite{Borisov} constructs the surface $Y$ as a system of quadrics in the variables $u_0,u_1, w_1, \ldots, w_6$ with the double fiber given by $\{u_1=0\}$. We compute a number of random points on the double fiber of $Y$ and use the equations to construct points on $\P^2_{Keum}$ lying on the preimage of the double fiber. We then look for polynomials vanishing on these points to compute $H^0(\P^2_{Keum}, 9H)$. The search for cubic polynomials gave 112 cubics with 16 in each $C_7$ weight.

\subsection*{Step 2: Constructing $|4H|$} \label{section2step2}
We may attempt to construct $|4H|$ as follows. The action of the $C_7$ automorphism subgroup on $H^0(\P^2_{Keum}, 4H)$ gives a $C_7$-representation which splits $H^0(\P^2_{Keum}, 4H)$ into three one-dimensional $C_7$-eigenspaces. The Holomorphic Lefschetz Fixed-Point formula shows that the eigenvalues are $\xi^3, \xi^5, \xi^6$, where $\xi$ is a seventh root of unity. Thus $H^0(\P^2_{Keum}, 4H)\cong \C r_3 \, \oplus\, \C r_5\, \oplus\, \C r_6 $ where $r_3, r_5, r_6$ are sections of $4H$ with $C_7$-weights 3, 5, and 6 respectively. In addition, the $C_3$-action on the surface implies $r_5=\sigma(r_3), r_6=\sigma^2(r_3)$ for $\sigma$ an order 3 automorphism on $\P^2_{Keum}$. The product $d=r_3 r_5 r_6$ is therefore a $C_3$-invariant section with $C_7$-weight 0 in $H^0(\P^2_{Keum}, 12H)$ (it is then invariant under the whole automorphism group). 

Set $s_i=r_i^3 \in H^0(\P^2_{Keum}, 12H)$ for $i\in \{3,5,6\}$. The equation
\[
s_3 s_5 s_6 = d^3
\]
 in $H^0(\P^2_{Keum}, 36H)$ allows us to narrow down parameters in the search for $r_3, r_5, r_6$. Since $s_3, s_5, s_6,$ and $d$ lie in $H^0(\P^2_{Keum}, 12H)$, they are quadratic in the variables $Q_0, \ldots, Q_9$ for the fake projective plane. It is sufficient to construct $s_3$ since $s_5$ and $s_6$ may be constructed from $s_3$ with the $C_3$ action. Additionally, since $s_3$ has $C_7$ weight $3 \times 3\equiv 2 \mod{7}$, we narrow the search down to $C_7$-weight 2 quadratics. 
 
We may further reduce the number of parameters with additional data. The curve $\{r_3=0\}$ passes through the two $C_7$ fixed points 
\begin{align*}
p_1&=(0\colon0\colon0\colon0\colon0\colon0\colon0\colon1\colon0\colon0),\\
p_2&=(0\colon0\colon0\colon0\colon0\colon0\colon0\colon0\colon1\colon0).\end{align*}
It follows that at these points the curve $\{s_3=0\}$ vanishes with multiplicity 3, which place additional conditions on the coefficients of $s_3$. 

Now we begin describing the details of the calculation. We first calculate the order 3 neighborhoods of the points $p_1$ and $p_2$. This was done by computing the tangent space and solving for the conditions of the neighborhoods vanishing on the FPP. We began by solving for the order 2 neighborhood and then for the third order. To speed up calculations, it was sufficient to take some equations for $\P^2_{Keum}$ locally cutting out the point. After computing these neighborhoods, we posit the general form for $s_3$ as weight 2 quadratics in the variables and then solve for the conditions of being identically 0 at the higher order neighborhoods. We are able to solve for two of these variables, narrowing down the general form for $s_3$ to 6 variables. 

We now want to solve for the sextic equation $s_3s_5s_6 -d^3=0$. The requirement that $d$ be invariant under the full automorphism group forces it to be of the form
\[
e_1 Q_0^2 + e_2 (Q_1 Q_6 + Q_2 Q_4 + Q_3 Q_5) + e_3 (Q_1 Q_9 + Q_2 Q_7 + Q_3 Q_8)
\]
for undetermined coefficients $e_1, e_2, e_3$. We also obtain the general forms for $s_5$ and $s_6$ by applying the $C_3$ automorphism to $s_3$. To solve for the coefficients, we compute some points of $\P^2_{Keum}$ with high accuracy and substitute them into $s_3s_5s_6 -d^3=0$ to obtain a system of 24 cubics in 6 variables. We solve this system of equations by applying the trick of \cite{Borisov and Fatighenti}. The Hilbert polynomial of the system of equations modulo 37 with $\sqrt{-7} \equiv 17 \mod{37}$ is 3, which suggests that there are 3 solutions for this system. By applying successive linear conditions on the system and checking the Hilbert polynomial at each step, we are able to take linear cuts that drop the Hilbert polynomial eventually to 1. At some point there are 3 different choices for the linear cuts corresponding to our 3 solutions. We were able to lift these 3 solutions modulo $37^{200}$ and then use the lattice reduction algorithm to obtain the corresponding solutions over $\Q(\sqrt{-7})$. The three solutions differed by a cube root of unity. We selected the solution defined over the desired field of definition to proceed.

The solution for these coefficients allow us to fully determine $s_3,s_5,s_6$, and $d$. The equations for $s_3$ and $d$ are given below, with $s_5$ and $s_6$ found by applying the $C_3$ automorphism. Points on $\{r_3=0\}$ may then be calculated by solving for the simultaneous conditions $\{s_3=0, d=0\}$. These points were used later in the construction.
\begin{align*}
 s_3 &=  \frac{\left(-212275+26525 i \sqrt{7}\right) Q_0 Q_5}{2470336}+\frac{\left(22575+51275 i \sqrt{7}\right) Q_0
   Q_8}{1235168}\\
   &+\frac{\left(139475+17575 i \sqrt{7}\right) Q_1 Q_2}{9881344}+\frac{\left(196875-91425 i \sqrt{7}\right)
   Q_3 Q_4}{2470336}\\
   &+\frac{\left(-303625-270725 i \sqrt{7}\right) Q_3 Q_7}{4940672}+\frac{\left(139475+17575 i
   \sqrt{7}\right) Q_6^2}{1235168}\\
   &+\frac{\left(795725-287175 i \sqrt{7}\right) Q_6 Q_9}{4940672}+\frac{\left(-57575-549675 i
   \sqrt{7}\right) Q_9^2}{9881344}\\
   d&=\frac{25}{9881344} \Bigl(3407 \sqrt{-7} Q_0^2+17045 Q_0^2-2812 \sqrt{-7} Q_1 Q_6-22316 Q_1 Q_6\\
   &+329 \sqrt{-7}
   Q_1 Q_9-21987 Q_1 Q_9-2812 \sqrt{-7} Q_2 Q_4-22316 Q_2 Q_4+329 \sqrt{-7} Q_2
   Q_7\\
   &-21987 Q_2 Q_7-2812 \sqrt{-7} Q_3 Q_5-22316 Q_3 Q_5+329 \sqrt{-7} Q_3 Q_8-21987
   Q_3 Q_8 \Bigr)
\end{align*}

\subsection*{Step 3: The map $\P^2_{Keum} \to \P^5$}
With the computations of $9H$ and $4H$ we may now find $5H$. We look at suitable linear combinations of the 112 polynomials vanishing on $18H-9H=9H$ additionally vanishing on $4H$ to obtain $18H-9H-4H=5H$.

We first compute some random points on $4H$ by solving for $\{d=0, s_3=0\}$ on the FPP.  $5H$ is then found by looking for linear combinations of the cubics defining $9H$ for each weight that vanish on these points. To verify that they are in $|5H|$ we also check that they do not vanish on the whole fake projective plane.

The six resulting degree 3 polynomials give us the map $\P^2_{Keum} \to \P^5$. We calculate points in the image of this map in $\P^5$ and find 56 degree 6 polynomials in new variables $Z_1, \ldots, Z_6$ that vanish at these points. These give the desired embedding of the fake projective plane.

\begin{rem}
The $C_7$-weights on the variables $Z_1, Z_2,\ldots, Z_6$ are $1,2,\ldots, 6$. There is no weight 0 variable. The construction required that we shift the $C_7$ weights by $3$. This may be explained by viewing our construction of $H^0(\P^2_{Keum},5H)$ as given by an embedding 
\[
H^0(\P^2_{Keum},5H) \hookrightarrow H^0(\P^2_{Keum}, 18H)
\]
with the map given by tensoring with $s_3 \otimes f$ for $s_3\in H^0(\P^2_{Keum}, 4H)$ and $f\in H^0(\P^2_{Keum}, 9H)$. While $f$ has weight 0, $s_3$ has weight $3$ and therefore shifts the weights of $H^0(\P^2_{Keum}, 5H)$ by 3. 
\end{rem}

We take care to reconstruct the automorphism group. While the $C_7$-action is preserved under our construction, the non-$C_3$-invariance of $s_3$ introduces a scaling factor in the $C_3$ action. We fix the coefficients of this scaling factor and recompute the equations with the scaling to find a better basis for the action. As before, we take random linear combinations of the equations that span the space and take the simplest ones to further simplify the equations.

Finally, we use Magma to verify that the Hilbert polynomial is as expected. The verification process for the FPP is carried out as in \cite{Borisov} working modulo $p=1327$ with $\sqrt{-7}=103 \mod{1327}$. Thus we have constructed Keum's fake projective plane as a degree 25 surface in $\P^5$. 

\section{Future Directions} \label{section4}

One hopes to find a coordinate change to additionally simplify the 56 equations in $\P^5$.

A related construction of interest is that of Mumford's original fake projective plane \cite{Mumford}. This surface has not been explicitly constructed yet. It lies in the same class as $\P^2_{Keum}$ and two other fake projective planes. We are currently attempting to find this surface by computing a seven-to-one cover of $\P^2_{Keum}$, after which several cover relations may yield the surface and the two fake projective planes in the same class.

\section{Acknowledgements}
The authors thank the DIMACS REU program at Rutgers University for supporting this research project. This work was carried out while the second author was supported by NSF grant CCF-1852215.


\end{document}